\documentclass[12pt]{amsart}
\usepackage{amssymb}
\usepackage{mathrsfs}

\theoremstyle{plain}

\newtheorem{theorem}{Theorem}
\newtheorem{corollary}{Corollary}[theorem] 
\newtheorem{problem}{Problem}
\newtheorem{remark}{Remark}
\theoremstyle{definition}

\newcommand{\dist}{{\rm dist\,}}

\newcommand{\co}{{\rm conv\,}}

\newcommand{\NN}{\mathbb N}

\newif\ifComplain
\Complaintrue           
\def\complain#1{\ifComplain\ifhmode \newline\fi{\sf *** \ \ #1
\\}\fi}

\newif\ifmarglab
\makeatletter
\def\label#1{\@bsphack\ifmarglab\marginpar{LAB:#1}\fi\if@filesw {\let\thepage\relax
   \def\protect{\noexpand\noexpand\noexpand}%
   \edef\@tempa{\write\@auxout{\string
      \newlabel{#1}{{\@currentlabel}{\thepage}}}}%
   \expandafter}\@tempa
   \if@nobreak \ifvmode\nobreak\fi\fi\fi\@esphack}
\makeatother
\marglabfalse

\usepackage{graphicx}

\long\def\onefigure#1#2{
\begin{figure*}[tbp]
\begin{center}
#1
\end{center}
\caption{#2}
\end{figure*}
} 

\newcommand\newipefig[2]
{\onefigure{\includegraphics{#1.pdf}}{\label{#1} #2}}

\begin{document}

\title[Weak limits]{Weak limits of consecutive projections and of greedy steps}

\author[P. A. Borodin]{Petr A. Borodin}
\address{Department of Mechanics and
Mathematics, Moscow State University, Moscow 119991, Russia}
\address{and}
\address{Moscow Center for Fundamental and Applied Mathematics}
\email{pborodin@inbox.ru}
\author[E. Kopeck\'a]{ Eva Kopeck\'a}
\address{Department of Mathematics\\
   University of Innsbruck\\
 A-6020 Innsbruck, Austria}
 \address{and}
\address{Moscow Center for Fundamental and Applied Mathematics}
\email {eva.kopecka@uibk.ac.at}




\begin{abstract}
Let $H$ be a Hilbert space. We investigate the properties of  weak  limit points of iterates of random projections onto $K\geq 2$ closed convex sets in $H$ and the parallel properties of weak  limit points of residuals of random greedy approximation with respect to $K$ dictionaries. In the case of convex sets these properties imply weak convergence in all the cases known so far. In particular, we give a short proof of the theorem of  Amemiya and Ando   on weak convergence when the convex sets are subspaces.
The question of the weak convergence in general remains open.
\end{abstract}

\maketitle

\date{\today}

\section*{Introduction}

In what follows $H$ is a real Hilbert space
with scalar product $\langle\,\cdot\,,\cdot\,\rangle$ and norm $|\,\cdot\,|$.

Let $A_1,\dots,A_K$ be closed and convex sets in  $H$, $K\geq 2$, so that  $A_1\cap\dots\cap A_K=\{0\}$.  Let $P_i$ denote the metric projection onto $A_i$. Let $i(n)\in \{1,\dots, K\}$ be a fixed sequence containing each $k\in \{1,\dots, K\}$ infinitely often. For $x_0\in H$, we consider the sequence
\begin{equation}\label{pr}
x_n=P_{i(n)} x_{n-1}, \qquad n=1,2,\dots
\end{equation}
In the case when $A_i$ are closed subspaces of $H$ the convergence properties of the sequence $\{x_n\}$ are well understood. If the sequence of the indices $\{i_n\}$ is periodic then the sequence $\{x_n\}$ converges in norm \cite{N}, \cite{Ha}. The rate of the convergence depending on the position of the subspaces and  of the initial point is known \cite{BaGM1}, \cite{BaGM2}, \cite{BDH}, \cite{BK}. In this context an interplay with the convergence properties of the greedy algorithm was discovered recently \cite{BK}.
If no extra information about the indices, or about the position  of the  subspaces is known
 already for $K=3$  divergence might occur \cite{P}, \cite{KM}, \cite{KP}, \cite{K20}. The sequence $\{x_n\}$, however, always converges weakly to zero \cite{AA}.

In the lack of linearity, when the sets $A_i$ are just closed and convex, the situation is different.
Already for $K=2$ the sequence $\{x_n\}$ might diverge in norm,  although the sequence of indices is inevitably periodic \cite{H},\cite{K08},\cite{MR}.
Weak convergence is known only under additional conditions: when $K\leq 3$ \cite{DR}, or when the indices are periodic \cite{Bre}, or when the sets are ``somewhat symmetric"   \cite{DR}.

We denote by
  $W=W(x_0)$ the set of all partial weak limits of the sequence (\ref{pr}) and  face
\begin{problem}\label{problem1}
 Is it true that $W=\{0\}$?
\end{problem}
We investigate the structure of the set $W$ and give new short proofs of the weak convergence in all
of the cases mentioned above. In particular, we give a short proof of the theorem of Amemiya  and Ando on  weak convergence when the convex sets are subspaces.
The general case remains, however, open.

In the spirit of \cite{BK}, we establish an interplay with the weak convergence problem of greedy approximation with respect to $K$ dictionaries. The structural properties of the set of weak partial limits of this greedy approximation turn out to be the same. We have hit the same bounds of knowledge while seeking   weak convergence.

\section{Projections on convex sets}\label{sec:convex}

Let $A_1,\dots,A_K$ be closed and convex sets in $H$, $K\geq 2$ so that  $A_1\cap\dots\cap A_K=\{0\}$. Let $i(n)\in \{1,\dots, K\}$ be a fixed sequence containing each $k\in \{1,\dots, K\}$ infinitely often and let the sequence $\{x_n\}$ be defined  by (\ref{pr}). We assume $i(n)\not=i(n+1)$ without loss of generality.

We study the structure of the set $W=W(x_0)$ of all partial weak limits of the sequence $\{x_n\}$.

Since the nearest point projection onto a convex set  is a $1$-Lipschitz mapping,  the norms
  $|x_n|$   decrease and the set $W$ is always nonempty. We may assume that $|x_n|\searrow R>0$, as  $R=0$ implies convergence in norm and hence $W=\{0\}$.

For $w\in W$, we denote by $J(w)$ the maximal subset of $\{1,\dots, K\}$ such that $w\in A_{J(w)}$. Here we use the notation $A_J=\cap_{j\in J}A_j$.

Since $|x_n-x_{n-1}|^2\leq |x_{n-1}|^2-|x_n|^2$, for $x_n\in A_{i(n)}$ we have
\begin{equation}\label{dist}
\dist (x_n,A_{i(n\pm m)})\to 0 \qquad (n\to \infty)
\end{equation}
for any fixed $m$. Therefore $|J(w)|\geq 2$ for each $w\in W$, and $W$ is a weakly closed subset of $\cup_{|J|\geq 2}A_J \cap B(0,R)$, where $B(0,R)$ is the ball centered at 0 of radius R.
If $w\not=0$, then $|J(w)|<K$, since $\bigcap A_i=\{0\}$.
It also follows from (\ref{dist}) that in case $i(n)\equiv n ({\rm mod}\, K)$ of alternating projections we have $J(w)=\{1,\dots,K\}$ for each $w$, and hence $W=\{0\}$. In particular,  if we have just  two convex sets then   the sequence $\{x_n\}$ converges weakly.

Next we show that if $W$ contains an element of maximal norm, then  $W=\{0\}$.

\begin{theorem}\label{theorem1}
 For each $w\in W$, $w\not=0$, one can find another element $w'\in W$ with the following properties:
\begin{itemize}
\item[(i)] $|J(w')\setminus J(w)|\geq 1$;
\item[(ii)]  $|J(w')\cap J(w)|\geq 2$
\item[(iii)]  $|J(w')|\geq 3$;
\item[(iv)]   $\langle w'-w,a\rangle\geq 0$ for every $a\in A_{J(w)}$.
\end{itemize}
In particular $|w'|>|w|$ in view of (i), since
$\langle w'-w,w\rangle\geq 0$, hence also $|w'|^2\geq |w|^2+|w-w'|^2$.
\end{theorem}
\begin{proof}
(i) Let
$$
x_{n_k}\rightharpoonup w, \qquad i(n_k)\in J(w).
$$
Taking a subsequence of $k$'s if needed, we can choose $q\notin J(w)$ with the following property:
for any $k$ there is a number $m_k\in (n_k, n_{k+1})$ with $i(m_k)=q$, so that
for any $n\in [n_k, m_k)$ we have $i(n)\in J(w)$, and hence $i(n)\not= q$.
Again taking a subsequence of $k$'s if needed, we get $x_{m_k}\rightharpoonup w'$, and this is the definition of $w'$.
Clearly $w'\in A_q$, hence $J(w')\ni q $ and (i) holds.

(ii) The numbers $i(m_k-1)$ and $i(m_k-2)$ belong to $J(w)$ and are distinct.  We choose two different numbers $i,j\in J(w)$ so that $i(m_k-1)=i$ and $i(m_k-2)=j$ for infinitely many $k$'s.  In view of (\ref{dist}) this implies $w'\in A_i\cap A_j$, hence $i,j \in J(w')\cap J(w) $ and (ii) holds.

The property (iii) follows from (i) and (ii).

(iv) For any $a\in A_{J(w)}$, we have
$$
\langle w'-w,a\rangle = \lim_{k\to \infty} \langle x_{m_k}-x_{n_k},a\rangle
$$
$$
= \lim_{k\to \infty} \sum_{n=n_k+1}^{m_k}\langle x_{n}-x_{n-1},a\rangle= \lim_{k\to \infty} \sum_{n=n_k+1}^{m_k-1}\langle x_{n}-x_{n-1},a\rangle
$$
$$
= \lim_{k\to \infty}\frac{1}{2} \sum_{n=n_k+1}^{m_k-1}(| x_{n-1}-a|^2-|x_n-a|^2 + |x_n|^2 -|x_{n-1}|^2)
$$
$$
=\frac{1}{2} \lim_{k\to \infty} \left(\sum_{n=n_k+1}^{m_k-1}(| x_{n-1}-a|^2-|P_{i(n)}x_{n-1}-P_{i(n)}a|^2 ) + |x_{m_k-1}|^2 -|x_{n_k}|^2\right)\geq 0,
$$
since every term in the sum is non-negative and $\lim_{k\to \infty}|x_{m_k-1}|=\lim_{k\to \infty}|x_{n_k}|= R$.
\end{proof}

\begin{remark}\label{remark1}
 The inequality (iv) holds for $a\in A_{J(w,w')}$, where $J(w,w')=\{i(n): n\in [n_k,m_k-1], k=1,2,\dots\}$.
 Since $J(w,w')\subset J(w)$, $A_{J(w,w')}$ can be strictly larger than  $A_{J(w)}$.
\end{remark}

The following corollary is a special case of Theorem~2 of \cite{DR}; our proof is different.

\begin{corollary}\label{cor1.1}
If $K\leq 3$, then $W=\{0\}$.
\end{corollary}
\begin{proof}
The case  $K=2$ we have explained above Theorem~\ref{theorem1}.
Assume that $K=3$ and that there is $w\in W\setminus \{0\}$. By Theorem \ref{theorem1} there is $w'\in W$ with $|w'|>|w|$ and $J(w')=\{1,2,3\}$  Hence $w'=0$ which is a contradiction.
\end{proof}

Assume all the convex sets $A_i$ are cones. Assume, moreover, that the intersection of any triple
of these cones with the unit sphere has a positive distance to the intersection of any other triple.
Then $W=\{0\}$ according to the next corollary.

\begin{corollary}\label{cor1.2}
Suppose for every $r>0$ there exists $\delta(r)>0$ so that for any two different triples $\{i,j,k\}$
and $\{i,j,l\}$ and elements $u\in A_{\{i,j,k\}}\cap S(0,r)$, $v\in A_{\{i,j,l\}}\cap S(0,r)$ we have $|u-v|>\delta(r)$. Then $W=\{0\}$.
\end{corollary}
\begin{proof}
 Suppose $W\not=\{0\}$.  Using Theorem \ref{theorem1}, we construct   a sequence $w_n\in W$ so that $w_1\not=0$, $w_{n+1}=w_n'$, $|J(w_n)|\ge 3$, $J(w_n)\not= J(w_{n+1})$ and $|J(w_n)\cap J(w_{n+1})|\ge 2$ for each $n$. So we get $w_n\in  A_{\{i,j,k\}} $ and $w_{n+1}\in  A_{\{i,j,l\}} $ for some $i$, $j$ and $k\not= l$ depending on $n$.
Since the sequence $|w_n|$ is bounded and increasing,
let $r=\lim_{n\to \infty} |w_n|$. 
Hence, $u_n=rw_n/(2|w_n|)\in  A_{\{i,j,k\}}$ for all sufficiently large $n$.
Denoting $w_n=(1+t_n)u_n$, $t_n>0$, for those $n$ we have
$$
\begin{array}{rcl}
|w_{n+1}|^2&\ge & |w_n|^2+|w_{n+1}-w_n|^2 =\\
&=& |w_n|^2+|(1+t_{n+1})u_{n+1}-(1+t_n)u_n|^2 \\
&\ge & |w_n|^2+|u_{n+1}-u_n|^2+2\langle u_{n+1}-u_n, t_{n+1}u_{n+1}-t_nu_n \rangle\\
&\ge &|w_n|^2+|u_{n+1}-u_n|^2+ 2(r/2)^2(t_{n+1}+t_n-(t_{n+1}+t_n))\\
&=&|w_n|^2+|u_{n+1}-u_n|^2>|w_n|^2+\delta(r/2)^2.
\end{array}
$$
That means, however,  that $|w_n|$ is unbounded.
\end{proof}

\begin{theorem}\label{theorem2}
 If $W\not=\{0\}$, then one can  find two different elements $w,w'\in W$ so that
$$
w={\rm weak}\lim_{k\to \infty} x_{n_k}, \qquad w'={\rm weak}\lim_{k\to \infty} x_{m_k},
$$
where $n_1<m_1<n_2<m_2<\dots$, and $i(n)\in J(w)\cap J(w')$ for any $n\in \cup_k (n_k,m_k)$. Consequently,

 $(i)$ $\langle w'-w,a\rangle\geq 0$ for every $a\in A_{J(w)}$,

$(ii)$ $\langle w'-w,b\rangle\geq 0$ for every $b\in A_{J(w')}$.

\end{theorem}

\begin{proof}
Both inequalities (i) and (ii) follow from the first statement of the Theorem.  The proof follows that of (iv) of Theorem \ref{theorem1}: all projections between $n_k$ and $m_k$ have indices from $J(w)\cap J(w')$.

To prove the first statement we take $w$ and $w'$ from Theorem \ref{theorem1}. All indices in $J=\{i(n): n\in \cup_k (n_k,m_k)\}$ belong to $J(w)$ by the proof of Theorem \ref{theorem1}. If $J\subset J(w')$, we are done. Otherwise we  define $\nu_k$ as the largest numbers  $n\in (n_k, m_k)$ such that $i(n)\notin J(w')$.
By taking a subsequence of $k$'s so that all these $i(n)$ are  the same, we get $x_{\nu_k}\rightharpoonup v\not= w'$. Then we redefine $w:=v$, $n_k:=\nu_k$. The renewed set $J=\{i(n): n\in \cup_k (n_k,m_k)\}$ is now  a subset of $J(w')$, and the number of elements in it has decreased by at least one. If this new $J$ is also included in the new $J(w)$, we stop. Otherwise we this time choose the numbers $\nu_k$ as the least numbers  $n\in (n_k, m_k)$ such that $i(n)\notin J(w)$. Then we redefine $w'$ and $m_k$'s. Since $|J|$ is decreasing, this oscillation process stops in a finite number of steps: $|J|$ cannot become less than 2. In case $|J|=2$  obviously   $J\subset J(w)\cap J(w')$.
\end{proof}

Dye and Reich used  in \cite{DR} the  so-called weak internal points (WIP) of a convex set to prove a
 nonlinear result that properly contains the
original linear theorem of Amemiya and Ando: if all the $K$ closed convex sets  are linear subspaces then the sequence $\{x_n\}$ converges weakly \cite{AA}. In  our version of the theorem we assume that
zero is a  WIP in   each of the  convex sets $A_k$. Again, the result is a special case of Theorem~5 of \cite{DR}; our proof is different.

\begin{corollary}\label{cor2.1}
Assume that zero is a weak internal point of each of the $K$ convex sets $A_k$: if $a\in A_k$ then
$-\lambda a\in A_k$ for some $\lambda=\lambda(a,k)>0$. Then $W=\{0\}$.
In particular, if all $A_k$ are closed linear subspaces of $H$ then the sequence (\ref{pr}) converges weakly.
\end{corollary}
\begin{proof}
Assuming $W\neq \{0\}$
we take the two different elements $w,w'\in W$ from Theorem \ref{theorem2}.
Using (i) of Theorem \ref{theorem2} for $a=w$ gives $\langle w'-w,w\rangle\geq 0$ and $|w'|^2\geq |w|^2+|w-w'|^2$.
Using (ii) of Theorem \ref{theorem2} for $b=-\lambda w'$ gives $\langle w'-w,-\lambda w'\rangle\geq 0$ and $|w|^2\geq |w'|^2+|w-w'|^2$.
Hence $w=w'$, which is a contradiction.

\end{proof}

\section{Parallels between projecting onto convex sets and  greedy approximation}

 A subset $D$ of the  the unit sphere $S(H)$ of the Hilbert space $H$ is called a {\em dictionary} if its span is dense in $H$.
 Assume, moreover, that $D$  does not lie in a half-space: for any nonzero $v\in H$, there exists $g\in D$ such that $\langle v,g\rangle>0$.
The greedy approximation algorithm then generates for $D$ and  for any element $x=x_0\in
H$  the sequence
\begin{equation}
\label{eq1}
x_{n+1}=x_n-\langle x_n,g_{n+1}\rangle g_{n+1}, \qquad n=0,1,\dots,
\end{equation}
where the element $g_{n+1}\in D$ is such that
$$
\langle x_n,g_{n+1}\rangle=\max\{\langle x_n,g\rangle\colon g\in D\}.
$$
The existence of $\max\{\langle x,g\rangle\colon g\in D\}$ for every $x\in H$ is an additional condition on $D$. If  the maximum is attained  on several elements of $D$, any of them is selected as $g_{n+1}$.
More precisely, this algorithm is called the \textit {pure greedy} algorithm, in contrast to other approximation algorithms whose names contain the word ``greedy'', see~\cite{T}.

For any symmetric dictionary $D$ the pure greedy algorithm converges in norm, see~\cite[Ch.~2]{T}.
That is, $x_n\to 0$ for any initial element $x=x_0$, and $x$ is represented as a norm-convergent series
$\sum_{n=0}^\infty \langle x_n,g_{n+1}\rangle g_{n+1}$. If $D$ is not symmetric, the greedy algorithm may diverge in norm~\cite{B21}, although it always converges weakly to zero~\cite{B20}.

Several details of the divergence construction in~\cite{B21} occur  to be   similar to that of~\cite{K08}.     The ``bridge" between this two seemingly different examples is the theorem of Moreau \cite{M}:
\begin{equation}\label{moreau}
P_A(x)=x-P_{A^*}(x)
\end{equation}
for any $x\in H$, any convex cone $A\subset H$ and its polar cone
$$
A^*=\{y\in H: \langle y, z\rangle\le 0 \quad \forall z\in A\}.
$$
Recall that both papers~\cite{H} and \cite{K08} provide examples of convex cones $A_1,A_2\subset H$ so that $A_1\cap A_2=\{0\}$ and alternating projections on those cones diverge in norm for certain starting elements. The formula (\ref{moreau}) allows us to interpret this result as an example of a divergent greedy algorithm with respect to the dictionary $D=(A_1^*\cup A_2^*)\cap S(H)$. Indeed, $D$ does not lie in a half-space as $A_1\cap A_2=\{0\}$, and for any greedy residual $x_n$ lying in, say, $A_1$, we have
$$
\max\{\langle x_n,g\rangle\colon g\in D\}=\max\{\langle x_n,g\rangle\colon g\in A_2^*\cap S(H)\},
$$
so that $x_{n+1}=x_n-P_{A_2^*}(x_n)=P_{A_2}(x_n)\in A_2$. Thus the author of~\cite{B21} didn't have to reinvent the wheel:~\cite{H} and \cite{K08} both provided the example he needed. 
However, the example in \cite{B21} is simpler than those of \cite{H} and \cite{K08}: it uses a discrete dictionary without the extra  care needed to build it of convex cones.

The above parallels between projecting onto convex sets and  greedy approximation have already been  noticed in \cite{BK} in the special case of subspaces. In the context of this paper, these parallels
 bring up  the question   of weak divergence of random greedy steps with respect to several dictionaries. This problem is considered in the next section. It turns  to have the same ``bounds of knowledge" as the problem of the weak divergence of random projections onto several  convex sets.

\section{Greedy approximation with respect to several dictionaries}

Let $K\ge 2$, $D_1,\dots,D_K$ be  subsets of $S(H)$ so that
  their union $\bigcup_{i=1}^K D_i$ is contained in no half-space: for any nonzero $v\in H$, there exists $g\in  \bigcup_{i=1}^K D_i$ such that $\langle v,g\rangle>0$. This implies that the set $\bigcup_{i=1}^K D_i$ is spanning; we will call here the sets $D_i$ dictionaries.

Assume that for each $x\in H$ and each $i\in \{1,\dots,K\}$ the following condition holds: if $\sup_{g\in D_i} \langle x, g \rangle> 0$, then the supremum  is attained on some element $g_i(x)\in D_i$. If it is attained at several elements of $D_i$, then
  we denote by $g_i(x)$ any one of them. If  $\sup_{g\in D_i} \langle x, g \rangle \le 0$, we put $g_i(x)=0$.

Clearly, our assumption means that the set $\Lambda(D_i)=\{\lambda g: \lambda\ge 0, g\in D_i\} $ is proximal, and the element $\langle x, g_i(x)\rangle g_i(x)$ belongs to the metric projection $P_{\Lambda(D_i)}(x)$.

Let $G_i$ denote the mapping corresponding to one step of the greedy algorithm with respect to the dictionary $D_i$:
$$
G_i(x)=x-\langle x, g_i(x)\rangle g_i(x).
$$
Note that
\begin{equation}\label{Pyth1}
|G_i(x)|^2= |x|^2- |x-G_i(x)|^2.
\end{equation}

Let $i(n)\in \{1,\dots, K\}$ be a fixed sequence containing each $k\in \{1,\dots, K\}$ infinitely often and such that $i(n)\neq i(n+1)$ for all $n\in \NN$. For $x_0\in H$, we consider the sequence
$$
x_n=G_{i(n)} x_{n-1}, \qquad n=1,2,\dots.
$$
As we have already mentioned above, this sequence may diverge in norm even in case of one dictionary. Both examples in \cite{H} and \cite{K08} can be interpreted as norm divergence examples of residuals $x_n$ for alternating greedy steps with respect to two dictionaries. So we are interested in weak convergence, just as in case of projections.
Denoting $W=W(x_0)$ the set of all partial weak limits of the sequence $\{x_n\}$, we face

\begin{problem}\label{problem2}
 Is it true that $W=\{0\}$?
\end{problem}

We may assume $x_n\not= x_{n-1}$ for all $n$, that is, $\sup_{g\in D_{i(n)}} \langle x, g \rangle> 0$.
According to (\ref{Pyth1}),
\begin{equation}\label{pyth}
|x_{n+1}|^2=|x_{n}|^2-|x_n-x_{n+1}|^2,
 \end{equation}
hence the norms $|x_n|$ are   decreasing. We may assume that $|x_n|\searrow R>0$, since  $R=0$ implies $W=\{0\}$.

We define the closed convex cones
$$
A_i=\{y\in H:  \langle y, g \rangle\le 0 \mbox{ for all } g\in D_i\}, \qquad i=1,\dots, K.
$$
Notice, that $A_i$ is the polar cone of $\overline{\co}\Lambda(D_i)$.
As in Section~\ref{sec:convex}, for $w\in W$, we denote by $J(w)$ the maximal subset of $\{1,\dots, K\}$ such that $w\in A_{J(w)}$, and again    use the notation $A_J=\cap_{j\in J}A_j$. Let us nevertheless  stress, that
the set $W$ is here the result of greedy approximation with respect to the dictionaries $D_1, \dots, D_K$.

Let us prove that $|J(w)|\geq 2$ for each $w\in W$. The convergence $x_{n_j}\rightharpoonup w$ implies the convergence $x_{n_j+m}\rightharpoonup w$
for any fixed $m$, since $\lim_{i\to \infty}|x_i-x_{i+m}|=0$ by  (\ref{pyth}). Suppose the sequence $i(n_j+m)$ contains some $k$ infinitely often. If $w\notin A_k$, then $ \langle w, g \rangle > \delta> 0$ for some $g\in D_k$, which yields $ \langle x_{n_j+m}, g \rangle > \delta$ for all sufficiently large $j$, so that $|x_{n_j+m+1}|^2\leq |x_{n_j+m}|^2- \delta^2$ for such $j$ with $i(n_j+m)=k$, and a contradiction with $|x_n|\searrow R>0$. So we get $w\in A_k$, and since one can find at least two such $k$'s using different $m$'s, we arrive at $|J(w)|\geq 2$.

The same argument shows that in case $i(n)\equiv n ({\rm mod}\, K)$ of alternating greedy algorithm we have $J(w)=\{1,\dots,K\}$ for each $w$, and hence $W=\{0\}$.

Thus, $W$ is a weakly closed subset of $\cup_{2\le |J|}A_J \cap B(0,R)$.

\begin{theorem}\label{theorem3}
 For each $w\in W$, $w\not=0$, one can find another element $w'\in W$ with the following properties:
\begin{itemize}
\item[(i)] $|J(w')\setminus J(w)|\geq 1$;
\item[(ii)]  $|J(w')\cap J(w)|\geq 2$
\item[(iii)]  $|J(w')|\geq 3$;
\item[(iv)]   $\langle w'-w,a\rangle\geq 0$ for every $a\in A_{J(w)}$.
\end{itemize}
In particular $|w'|>|w|$ in view of (i), since
$\langle w'-w,w\rangle\geq 0$, hence also $|w'|^2\geq |w|^2+|w-w'|^2$.
\end{theorem}

\begin{proof}
 Theorem \ref{theorem3} is formally identical to Theorem \ref{theorem1}, and the proofs of (i)-(iii)  follow the same reasoning.

The proof of (iv) is slightly different. As in the proof of Theorem \ref{theorem1}, we have
two alternating sequences $n_1<m_1<n_2<m_2<\dots$ so that
$$
x_{n_k}\rightharpoonup w, \qquad x_{m_k}\rightharpoonup w',
$$
and $i(n)\in J(w)$ for all $n\in \cup_k[n_k,m_k)$.

For any $a\in A_{J(w)}$, we have
\begin{equation}\notag
\begin{split}
\langle w'-w,a\rangle &= \lim_{k\to \infty} \langle x_{m_k}-x_{n_k},a\rangle \\
&= \lim_{k\to \infty} \sum_{n=n_k+1}^{m_k}\langle x_{n}-x_{n-1},a\rangle= \lim_{k\to \infty} \sum_{n=n_k+1}^{m_k-1}\langle x_{n}-x_{n-1},a\rangle \\
 &= \lim_{k\to \infty}\sum_{n=n_k+1}^{m_k-1} (-1)\langle x_{n-1}, g_{i(n)}(x_{n-1})\rangle \langle g_{i(n)}(x_{n-1}),a\rangle
\geq 0.
\end{split}
\end{equation}
The last inequality holds since each of the summands is non-negative:
$\langle x, g_i(x)\rangle\ge 0$ for any $x$ and $i$ by the definition of $g_i$,  and  $\langle g_{i(n)}(x_{n-1}),a\rangle\le 0$  since  $i(n)\in J(w)$ and $a\in A_{J(w)}$.
\end{proof}

\begin{remark}\label{remark2}
 The inequality (iv) holds for $a\in A_{J(w,w')}$, where $J(w,w')=\{i(n): n\in [n_k,m_k-1], k=1,2,\dots\}$.
 Since $J(w,w')\subset J(w)$, $A_{J(w,w')}$ can be strictly larger than  $A_{J(w)}$.
\end{remark}

\begin{corollary}\label{corol3.1.}
If $K\leq 3$, then $W=\{0\}$.
\end{corollary}
\begin{proof}
If $K=2$ we have an alternating greedy algorithm,
hence convergence as we have explained above Theorem~\ref{theorem3}.

Assume that $K=3$ and that there is $w\in W\setminus \{0\}$. By Theorem \ref{theorem3} there is $w'\in W$ with $|w'|>|w|$ and $J(w')=\{1,2,3\}$  Hence $w'=0$ which is a contradiction.
\end{proof}

\begin{corollary}\label{cor3.2}
Suppose for any four indices  $i,j,k,l\in \{1,\dots,K\}$ the inequality
 \begin{equation}\label{ijkl}
 \inf_{s\in S(H)}\sup_{g\in D_i\cup D_j\cup D_k\cup D_l} \langle s, g\rangle>0
 \end{equation}
 holds. Then $W=\{0\}$.
\end{corollary}
\begin{proof}
The inequalities (\ref{ijkl}) provide $\delta>0$ so that for any distinct $i,j,k,l$ and $u\in A_{\{i,j,k\}}\cap S(H)$ there exists $g\in D_l$ such that $\langle u,g \rangle>\delta$. Hence, for any two different triples $\{i,j,k\}$
and $\{i,j,l\}$ and unit elements $u\in A_{\{i,j,k\}}$, $v\in A_{\{i,j,l\}}$ we have $|u-v|>\delta$:
$$
|u-v|\ge \langle u-v, g \rangle \ge \langle u,g \rangle>\delta.
$$

Further we repeat the proof of Corollary~\ref{cor1.2}.
Suppose $W\not=\{0\}$.  By Theorem \ref{theorem3}, we can produce a sequence $w_n\in W$ so that $w_{n+1}=w_n'$, $|J(w_n)|\ge 3$, $J(w_n)\not= J(w_{n+1})$ and $|J(w_n)\cap J(w_{n+1})|\ge 2$ for each $n$. So we get $w_n\in  A_{\{i,j,k\}} $ and $w_n\in  A_{\{i,j,l\}} $ for some $i$, $j$ and $k\not= l$ depending on $n$. Therefore, using that the sets $A$ are cones, we can refine the inequality from Theorem \ref{theorem3}:
$$
|w_{n+1}|^2\ge |w_n|^2+|w_{n+1}-w_n|^2\ge |w_n|^2+\delta^2|w_n|^2/2.
$$
That, however, means that $|w_n|$ is unbounded.

\end{proof}

\begin{theorem}\label{theorem4}
If $W\not=\{0\}$, then one can  find two different elements $w,w'\in W$ so that
$$
w={\rm weak}\lim_{k\to \infty} x_{n_k}, \qquad w'={\rm weak}\lim_{k\to \infty} x_{m_k},
$$
where $n_1<m_1<n_2<m_2<\dots$, and $i(n)\in J(w)\cap J(w')$ for any $n\in \cup_k (n_k,m_k)$. Consequently,

 $(i)$ $\langle w'-w,a\rangle\geq 0$ for every $a\in A_{J(w)}$,

$(ii)$ $\langle w'-w,b\rangle\geq 0$ for every $b\in A_{J(w')}$.

\end{theorem}

\begin{proof}

We repeat the proof of Theorem~\ref{theorem2}; it is purely combinatorial. The inequalities follow from the first statement  as in the proof of part (iv) of Theorem~\ref{theorem3}.
\end{proof}

\begin{corollary}\label{cor4.1}
If all $D_k$ are symmetric, then $W=\{0\}$.
\end{corollary}
\begin{proof}

Assume that $W\neq \{0\}$.
We take the two different elements $w,w'\in W$ from Theorem~\ref{theorem4}.
 Using (i) of Theorem~\ref{theorem4} for $a=w$ gives $\langle w'-w,w\rangle\geq 0$, hence $|w'|^2\geq |w|^2+|w-w'|^2$.
Using (ii) of Theorem~\ref{theorem4} for $b=-\lambda w'$ gives $\langle w'-w,-\lambda w'\rangle\geq 0$, hence $|w|^2\geq |w'|^2+|w-w'|^2$.
Thus we get $w=w'$, which is a contradiction.
\end{proof}

\end{document}